\documentclass[12pt,leqno]{article}
\usepackage{amsmath}
\usepackage{amssymb}
\usepackage{theorem}
\usepackage{epsf}
\usepackage[scanall]{psfrag}
\usepackage{graphicx}

\numberwithin{equation}{section}

\theoremstyle{plain}
\newtheorem{thm}{Theorem}[section]
\newtheorem{cor}[thm]{Corollary}
\newtheorem{lem}[thm]{Lemma}

\makeatletter
\DeclareRobustCommand{\qed}{%
  \ifmmode 
  \else
\leavevmode\unskip\penalty9999 \hbox{}\nobreak\hfill
  \fi
  \quad\hbox{\qedsymbol}}
\newcommand{\openbox}{\leavevmode
  \hbox to.77778em{%
  \hfil\vrule
  \vbox to.675em{\hrule width.6em\vfil\hrule}%
  \vrule\hfil}}
\newcommand{\qedsymbol}{\openbox}
\newenvironment{proof}[1][\proofname]{\par
  \normalfont
  \topsep6\p@\@plus6\p@ \trivlist
  \item[\hskip\labelsep\itshape
    #1\@addpunct{.}]\ignorespaces
}{%
  \qed\endtrivlist
}
\newcommand{\proofname}{Proof}
\newcommand{\references}[1]{\bigskip\bigskip\begin{center}REFERENCES\end{center}
{#1}\bigskip}
\newcommand{\Ref}[1]{\par{\vbox{\hangindent=3pc \hangafter=1 {\noindent {#1}} }\smallskip}}

\makeatother

\def\D{\displaystyle}

\font\tenmsbm=msbm10
\font\sevenmsbm=msbm7
\font\fivemsbm=msbm5
\newfam\msbmfam \textfont\msbmfam=\tenmsbm \scriptfont\msbmfam=\sevenmsbm
 \scriptscriptfont\msbmfam=\fivemsbm

\def\hexnumber@#1{\ifcase#1 0\or 1\or 2\or 3\or 4\or 5\or 6\or 7\or 8\or
 9\or A\or B\or C\or D\or E\or F\fi}

\newcommand{\z}{\mathbb{Z}}
\newcommand{\bZ}{\mathbb{Z}}
\newcommand{\bG}{\mathbb{G}}
\newcommand{\zm}{\z_m}
\newcommand{\walk}{G(m, p, p')}
\newcommand{\ggwalk}{G(m, \bp,\bq)}
\newcommand{\gwalk}{G(m, \bp)}
\newcommand{\e}{\varepsilon}
\newcommand{\boldu}{\mathbf{u}}
\newcommand{\bolds}{\mathbf{s}}
\newcommand{\boldt}{\mathbf{t}}
\newcommand{\boldv}{\mathbf{v}}
\newcommand{\boldp}{\mathbf{p}}
\newcommand{\boldq}{\mathbf{q}}

\newcommand{\boldone}{\mathbf{1}}

\def\sqr#1#2{{\vcenter{\vbox{\hrule height.#2pt
        \hbox{\vrule width.#2pt height#1pt \kern#2pt
        \vrule width.#2pt}
        \hrule height.#2pt}}}}

\def\sqr#1#2{{\vcenter{\vbox{\hrule height .#2pt
        \hbox{\vrule width.#2pt height#1pt \kern#1pt
        \vrule width.#2pt}
        \hrule height.#2pt}}}}
\def\squar4{\mathchoice\sqr74\sqr74\sqr{2.1}3\sqr{1.5}3}
\def\={{{{\mathop{=}\limits^{L}}}}}

\def\e{{{\varepsilon}}}

\def\bp{\mathbf{p}}
\def\bq{\mathbf{q}}
\def\bP{\mathbf{P}}
\def\bQ{\mathbf{Q}}

\def\PP{\mathbb{P}}

\def\bC{\mathbb{C}}
\def\bH{\mathbb{H}}
\def\bK{\mathbb{K}}
\def\bI{\mathbb{I}}
\def\bfalpha{\boldsymbol{\alpha}}
\def\bfbeta{\boldsymbol{\beta}}
\def\bdelta{\boldsymbol{\delta}}
\def\brho{\boldsymbol{\rho}}
\def\bmu{\boldsymbol{\mu}}
\def\btau{\boldsymbol{\tau}}
\def\bpi{\boldsymbol{\pi}}
\def\boldb{\mathbf{b}}
\def\bfx{\mathbf{x}}

\begin{document}

\thispagestyle{empty}

\begin{titlepage}
\begin{center}
{\Large On Random Walks and Diffusions Related to Parrondo's Games}
\end{center}
\vskip.5in
\begin{center}
{\large Ronald Pyke}\\
{\large University of Washington}\\
{\large\quad}\\
{\large\quad}
\end{center}
\vskip.5in
\begin{center}
\today
\end{center}
\vskip.5in
\abstract{
In a series of papers, G. Harmer and D. Abbott
study the behavior of random walks associated with games introduced in 1997
by J. M. R. Parrondo. These games illustrate an apparent paradox that
random and deterministic mixtures of losing games may produce winning games.
In this paper, classical cyclic random walks on the
additive group of integers modulo $m$, a given integer, are used
in a straightforward
way to derive the strong law limits of a general class of games that contains
the Parrondo games. We then consider the question of when random mixtures
of fair games related to these walks may result in winning games.
Although the context for these problems is elementary, there remain open
questions. An extension of the structure of these walks to a class
of shift diffusions is also presented, leading to the fact that a random
mixture of two fair shift diffusions may be transient to $+\infty$. }

\vfill


\end{titlepage}

\bigskip

\section{Introduction}
The purpose of this paper is to study a family of random walks that include
those arising in the games devised by J. M. R. Parrondo
in 1997 to illustrate the
apparent paradox that two `losing' games can result in a `winning' game when
one alternates between them.
We refer the reader to Harmer and Abbott (1999a,b),
Harmer, Abbott and Taylor(2000) and Harmer, Abbott, Taylor and Parrondo(2000)
in which Parrondo's paradox is discussed, large simulations of specific Parrondo
games and mixtures thereof are presented and certain theoretical
results are given. These authors also give a
heuristic explanation of the paradox in terms of the Brownian ratchet,
the original motivation for the suggestion of these games. Other references
to the general subject are included in the above mentioned papers by
Harmer and Abbott. The reader may also note the reference Durrett, Kesten
and Lawler(1991) which also deals with the general question of showing
that winning games can be formed by mixing fair ones.

The suggested paradox may be visualized as follows. You are about to
play a two-armed slot machine. The casino that owns this two-armed bandit
advertises that both arms on their two-armed machines are "fair" in the sense
that any player who plays either of the arms is assured that the average
cost per play approaches zero as the number of plays increase. However,
the casino does not constrain you to stay with one arm; you are allowed
to use either arm on every play. You just tell the machine before
beginning how many plays you wish to make. At the
end of that number of plays, the machine displays the total amount won or
lost. The question of interest in this context would be whether it is
possible for the casino to still make money using only "fair" games.

In this paper a random walk will refer to a Markov chain $\{S_n:n=0,1,2,
\ldots\}$ taking values in the integers, $\z$, which satisfies the discrete
continuity condition
\[
|S_n-S_{n-1}|=1\quad\hbox{a.s.\ for each}\quad n\ge 1.
\]
Let the transition probabilities for the random walk be denoted by
\[
p_j=P(S_{n+1}-S_n=1|S_n=j),\qquad q_j=P(S_{n+1}-S_n=-1|S_n=j)
\]
and
\[
r_j=1-p_j-q_j=P(S_{n+1}=S_n|S_n=j)
\]
for $j\in\bZ$. Assume that $p_jq_j\ne 0$ for all $j$.
For fixed integer $m\ge 1$, define a \emph{mod m random walk} to be a
random walk in which the transition probabilities $p_i$, $r_i$, $q_i$ depend
only upon the congruence class  $\mod{m}$ of the state $i$.
Thus, these {\it lattice regular} or {\it periodic} random walks are such
that for some specified integer $m>1$, $p_j=p_{j+m}$ and $q_j=q_{j+m}$
for all $j\in \z$.
More generally, define a \emph{mod m Markov chain} on the integers
to be one whose parameters
depend only upon the congruence classes  $\mod{m}$ of the states, namely,
$p_{ij} = p_{i+m,j}$ for all integers $i,j$.
This paper is concerned with the case of random walks,
but places where the approach applies more generally are pointed out.

A $\mod{m}$ random walk is determined by the $2m$ parameters
$p_j, q_j;0\le j <m$.
Write $\bp=(p_0,p_1,\cdot\cdot \cdot,p_{m-1})$
with an analogous use of $\bq$ to specify the walk's parameters. Observe that when $m=1$ the walk is
classical simple random walk, so our main interest is in the cases of $m>1$.

These random walks are viewed as games with the increment
$X_n=S_n-S_{n-1}(n\ge 1)$ denoting the
gain at the $n$-th play. We say that the game is a winning/losing/fair game
according as the almost sure limit of $S_n/n$ is positive/negative/zero.

For given $m$, write $\zm:=m\z=\{km:k\in\z\}$ for the integer lattice
of span $m$. In the games introduced by Parrondo, it is assumed
that the transition probabilities depend on the state only to the extent
that the state is or is not in $\zm$. Thus,
Parrondo's games are characterized by
\begin{equation}\label{eqn1.1}
P(X_{n+1}=1\mid S_0,S_1,\ldots, S_n) = p'1_{[S_n\in\zm]}+p1_{[S_n\not\in\zm]}
\end{equation}
for some $p, p'\in[0,1]$ and all $n\ge 0$. Write $q=1-p$ and $q'=1-p'$.
We may also write $k\equiv j\mod{m}$ when $k\in j+\zm$.

For simplicity, we write $\ggwalk$ to denote a general $\mod{m}$
random walk or game, but write $\gwalk$ for the game when each $q_j=1-p_j$
(i.e. each $r_j=0$) and write $\walk$ for the special Parrondo
random walk or game satisfying (\ref{eqn1.1}).

The required notation and preliminary structure are introduced in the
following section, in which the limiting results for $\walk$ games are
given for illustration. The general case is covered in Section 3, while in
Section 4 we resolve the central question about whether random mixtures
of losing Parrondo's games can be winning ones. The asymptotic gain is
derived in Section 5 while in Section 6 a certain expected interoccurrence time
that appears in the previously obtained expression for this is also dreived.
The method used to solve the recursion equations in these sections makes use
of an extension of results of Mihoc and Fr\'echet (cf. Fr\'echet(1952)) that are
provided in the Appendix to this paper. Continuous analogues to the
random walks considered here are introduced in Section 7. These $\mod{m}$
diffusions have drift functions that are periodic step functions so that
their embedded walks on the integers are $\ggwalk$ walks. In Theorem 7.2
the drift rates under which the embedded walk
has specified transition probabilities is determined.

\section{Preliminaries and Parrondo's Examples }

In the games suggested by Parrondo, the transition probabilities depend
on the state only  to the extent that it is or is not in $\zm$; see
(\ref{eqn1.1}) above.
The asymptotic behavior of these games, as for any $\mod{m}$
random walk  is determined by that of
its embedded walk on the lattice $\zm$. Since this embedded walk
is equivalent to simple random walk, its asymptotics are
well known and dependent solely upon a single parameter, the walk's
probability of 'success'. In this section we introduce the notation
required for the general case in Section 3 below, and illustrate the approach
in the special case of a Parrondo $\walk$ walk by substituting in known results
for simple random walk.

Let $T_1<T_2<\cdots$
be the successive transition times of the embedded walk on $\zm$. That is
$T_1 =\min\{n\ge 0:S_n\in \zm\}$ and, for $k>1$,
\[
T_{k+1}=\min\{n>T_k:S_n-S_{T_k}=\pm m\}
\]
with the minimum of a null set being defined to equal $+\infty$. Set
$T_0=0$. Write
\[
J_n=m^{-1}S_{T_{n+1}},\ \  \xi_n=T_{n+1}-T_n
\]
for $n\ge 0$ so that $\{(J_n, \xi_n):n\ge 0\}$ is a (possibly delayed)
Markov renewal process (MRP) in which the embedded random walk $\{J_n\}$
is simply a classical random walk
with constant probability of ` success',
\begin{equation}\label{eqn2.1}
p^\ast_m:=P(J_{n+1}-J_n=1\mid J_n).
\end{equation}
Hence, once
$p_m^{\ast}$ is known, the winning/losing/fair nature of the walk is
easily determined.

In general,
for $n$ satisfying $T_k<n\le T_{k+1}$,
\[
\frac{mJ_{k+1}-2m}{T_{k+1}} = \frac{S_{T_{k+1}}-2m}{T_{k+1}} \le
\frac{S_n}{n} \le \frac{S_{T_k}+m}{T_k} = \frac{mJ_k+m}{T_k}.
\]
Thus,
\begin{equation}\label{eqn2.7}
m\left(\frac{mJ_{k+1}}{k+1} - \frac{2}{k+1}\right)/\frac{T_{k+1}}{k+1} \le
\frac{S_n}{n} \le m\left(\frac{J_k}{k} + \frac{1}{k} \right)/\frac{T_k}{k}.
\end{equation}
It is known for the classical random walk $\{J_n\}$ that $J_k/k$ converges
a.s.\ as $k\to\infty$ to $p^\ast_m-q^\ast_m$, with $q^\ast_m = 1-p^\ast_m$.
Moreover, the stopping times
$\{T_k\}$ are partial sums of iid r.v.'s having finite expectations so that
\[
T_k/k \mathop{\longrightarrow}\limits^{\mathrm{a.s.}} E(T_2-T_1)<\infty.
\]
Upon taking limits in (\ref{eqn2.7}) one obtains
that with probability one,
\begin{equation}\label{eqn2.8}
\lim_{n\to\infty} \frac{S_n}{n} = \frac{m(p^\ast_m-q^\ast_m)}{E(T_2-T_1)}.
\end{equation}
Clearly then,  this limit is $0$, $>0$ or
$<0$ according as $p^\ast_m=$, $>$ or
$<q^\ast_m$.

The quantity $p^\ast_m$ is evaluated for the general $\ggwalk$ walk
in Lemma \ref{lem6.2} below.
However, for the special Parrondo $\walk$ random walk, the evaluation
is immediate once we introduce the notation and approach that is
needed for the general case, and
so we give it separately here as

\begin{lem}\label{lem2.1}
The $\zm$-embedded MRP of the $\walk$ random walk has transition probabilities
determined by the `success' probability
\begin{equation}\label{eqn2.2}
p^\ast_m=\frac{p'p^{m-1}}{p'p^{m-1}+q'q^{m-1}}
\end{equation}
for all $p,p'\in[0,1]$ satisfying $|p-p'| < 1$.
\end{lem}

\begin{proof}
The first part of this proof, through (\ref{eqn2.4}) below, is general and will
be needed in Section 3. The rest is substitution of known results.

Suppose $J_n=k$. That is, for the original walk suppose
$S_{T_{n+1}}=km$. Since $T_{n+1}$ is a stopping time, $P(J_{n+1}-J_n=1\mid J_n=k
)$
is just the probability that starting at $S_0=0$, the random walk $\{S_n\}$
reaches $m$ before it reaches $-m$. But $S_1$ equals 1 or $-1$ with
probability $p_0$ or $q_0$, respectively.
Thus if we let $A$ denote the event
that $\{S_n:n>1\}$ reaches 0 before it reaches $mS_1$ then the Markov
property implies that  $p^\ast_m$, the success probability
for the embedded walk, satisfies the following recursion relation,
in which we partition the event according to whether the original walk hits
zero before $m$ or not:
\begin{equation}\label{eqn2.3}
p^\ast_m=P(A)p^\ast_m+p_0\{1-P(A\mid S_1=1)\}.
\end{equation}
Hence
\begin{equation}\label{eqn2.4}
p^\ast_m=p_0P(A^c\mid S_1=1)/P(A^c).
\end{equation}
Since for the special case of this lemma, the conditional probabilities
given $S_1$ are just those that arise in
the classical gambler's ruin problem, (cf Feller (1968, Chap.\ XIV)
it is known that
\begin{equation}\label{eqn2.5}
P(A^c \mid S_1=1) = \left\{\begin{array}{ll}
(qp^{m-1}-p^m)/(q^m-p^m) & \hbox{if }p\ne q \\
1/m & \hbox{if }p=\frac{1}{2}\end{array}\right.
\end{equation}
and $P(A^c\mid S_1=-1)$ is similar but with $p$ and $q$ interchanged.
Substitution of (\ref{eqn2.5}) into (\ref{eqn2.4})  now
gives, when $p\ne q$,
\begin{equation}\label{eqn2.6}
P(A^c)
 =  (q'q^{m-1}+p'p^{m-1})(q-p)/(q^m-p^m)
\end{equation}
and, therefore, $p^\ast_m$  is as required by (\ref{eqn2.2}).
When $p=\frac{1}{2}$, the  substitution of (\ref{eqn2.5}) yields
$p^\ast_m=p'$ to complete the proof.
\end{proof}

Note that by (\ref{eqn2.2}), $p^\ast_m$
is the conditional probability that $S_n$ reaches $m$ before $-m$ given
that $S_0=0$ and that the first $m$ steps of $S_n$ are monotone.
This structure is more readily seen in the general case of Lemma \ref{lem6.2}
below.

For the special $\walk$ case the above result yields
\begin{cor}\label{cor2.2}
(Harmer and Abbott(2000a))
When $|p-p'| < 1$, the game $\walk$ is a fair, winning or losing game according
as
\[
p'p^{m-1}-q'q^{m-1}=0,\qquad >0\qquad \hbox{or}\qquad <0.
\]
\end{cor}

The condition in Corollary \ref{cor2.2} is more clearly expressed in terms
of new variables $x=p/q$ and $y=p'/q'$, namely, the game $\walk$ is a fair,
winning or losing one according as
\begin{equation}\label{eqn2.9}
y-x^{-(m-1)}=0,\qquad >0\qquad\hbox{or}\qquad <0.
\end{equation}
Recall that the degenerate case $q'=0=p$ has been excluded. Since the
inverse relationships are  $p=x/(1+x)$ and
$p'=y/(1+y)$, it follows from (\ref{eqn2.9}) that $\walk$ is fair if for
some $x\ge 0$, $p$ and $p'$ are related as
\[
q=1-p=\frac{1}{1+x}\qquad\hbox{ and }\qquad
p'=\frac{x^{-m+1}}{1+x^{-m+1}}=\frac{1}{1+x^{m-1}}.
\]

Here are some examples. For $x=1$, $G(m,\frac{1}{2}, \frac{1}{2})$ is fair
for every $m\ge 1$. For $m=4$ and $x=2$, the game $G(4, \frac{4}{5},
\frac{1}{65})$ is seen to be fair, and for $m=5$ and $x=2$, $G(5, \frac{2}{3},
\frac{1}{17})$ is fair. When $m=3$ and one chooses $x=3$, one
obtains the fair game $G(3, \frac{3}{4}, \frac{1}{10})$. The associated
games $G(3,
\frac{3}{4}-\e, \frac{1}{10}-\e)$ for a range of $\e>0$ are the losing
games used in the
simulation study of Harmer and Abbott(1999a).
The fact that these are losing games as indicated there is immediate from
the following observation: If $G(m, p_0, p'_0)$ is a fair game, then $\walk$
is a losing game whenever $0\le p'\le p'_0$ and $0\le p\le p_0$ with
 $p+p'<p_0+p'_0$; simply observe that $p'/q'$ and $p/q$ are increasing
functions of $p'$ and $p$, respectively, so that
$p'_0/q'_0=(p_0/q_0)^{-m+1}$ implies  $p'/q'\le (p/q)^{-m+1}$ whenever
$p'\le p'_0$ and $p\le p_0$. Since $G(3, \frac{3}{4}, \frac{1}{10})$ is a
fair game the result follows.

\section{General Mod m Random Walks}

Let $\{S_n:n\ge 0\}$ be a general (discretely continuous)
random walk on the integers
$\bZ$ in the sense described in the Introduction above.
The asymptotic
behavior of $\{S_n\}$ can be described in terms of the two associated
reflecting random walks on the negative and positive integers. The latter is
obtained, for example, by replacing $r_0$ and $q_0$ by $\bar r_0=1-p_0$ and
$\bar q_0=0$. It is known (cf.\ Feller (1968), Chap.\ XV.8 or Chung (1967),
Sect. I.12) that the
corresponding reflecting random walk on $\bZ^+=\{0,1,2,\ldots\}$ is
recurrent or transient according to
\begin{equation}\label{eqn6.1}
\sum^\infty_{i=1}\frac{q_1q_2\cdots q_i}{p_1p_2\cdots p_i}=\infty
\end{equation}
or not. When one looks similarly at the reflecting random walk on
$\bZ^-=\{0, -1, -2,\break
\ldots\}$, the roles of the $p$'s and $q$'s are
interchanged so that recurrence in this case holds if and only if
\begin{equation}\label{eqn6.2}
\sum^\infty_{i=1} \frac{p_{-1}p_{-2}\cdots p_{-i}}{q_{-1}q_{-2}\cdots
q_{-i}}=\infty.
\end{equation}
Now return to the original walk on $\bZ$. The positive part of this walk,
$\{S^+_n\}$, is a Markov renewal process in which all sojourn times are
equal to one except those between successive visits to state $0$. The
distribution of these latter sojourn times is a possibly deficient mixture that
includes with probability $q_0$ the distribution of the first passage time
from state $-1$ to state $0$. The latter passage time is finite with
probability one only if the reflecting random walk on $\bZ^-$ is recurrent.
Hence $\{S_n\}$ is a recurrent random walk if and only if both reflecting
random walks are recurrent, or equivalently, if and only if
both (\ref{eqn6.1}) and (\ref{eqn6.2}) hold. Consequently, the walk is
transient if and only if at least one of these series converges. Accordingly,
the boundary of a transient random walk may consist of either or both of
$+\infty$ and $-\infty$, depending upon which one or both of the series
converge. (Cf.\ Karlin and McGregor (1959), Section 4 where the integral
representations of the transition probabilities of the doubly infinite
random walk are expressed in terms of those of the two corresponding
reflecting walks.)

Consider now, for fixed integer $m\ge 1$, a \emph{mod m random walk} as
defined in Section 1.
random walk in which the transition probabilities $p_i$, $r_i$, $q_i$ depend
only upon the congruence class  $\mod{m}$ of the state $i$. (Note that when
$m=1$, the $\mod{1}$ random walk is just the classical random walk with
constant transition probabilities.)
 Thus, for $i=sm+l$ for some $s\in\bZ$
and $l=0,1,\ldots, m-1$, we know that
$(p_i, r_i, q_i)=(p_l, r_l, q_l)$. Moreover,
for $s\ge 0$, the summand in (\ref{eqn6.1}) becomes
\begin{equation}\label{eqn6.3}
\frac{q_1q_2\cdots q_i}{p_1p_2\cdots p_i} =
\frac{p_0}{q_0}\left\{\frac{q_0q_1\cdots q_{m-1}}{p_0p_1\cdots
p_{m-1}}\right\}^s \frac{q_0q_1\cdots q_l}{p_0p_1\cdots p_l}
\end{equation}
while for $s<0$, a similar representation holds with the $p$'s and $q$'s
interchanged. If we define
\begin{equation}\label{eqn6.4}
\rho_m:=\frac{p_0p_1\cdots p_{m-1}}{q_0q_1\cdots q_{m-1}}
\end{equation}
then the divergence of (\ref{eqn6.1}) holds if and only if $\rho_m\le 1$ while
(\ref{eqn6.2}) holds if and only if $\rho_m\ge 1$. By the above discussion,
the walk is then
recurrent, transient to $+\infty$ or transient to $-\infty$ according as
$\rho_m$ is equal to, greater than or less than one. This then proves

\begin{lem}\label{lem6.1}
For $m\ge 1$, a $\mod{m}$ random walk is recurrent, transient toward
$+\infty$ or transient toward $-\infty$ according as
\begin{equation}\label{eqn6.5}
p_0p_1\cdots p_{m-1}-q_0q_1\cdots q_{m-1}=0,\qquad >0 \quad \hbox{ or }<0.
\end{equation}
\end{lem}

It remains to evaluate $p^{\ast}_m$, the probability of 'success',
$p^{\ast}_m$, for the
embedded walk on $\zm$.

\begin{lem}\label{lem6.2}
For $m\ge 1$, and a $\mod{m}$ random walk $\ggwalk$ with parameters
\break $\boldp=(p_0,p_1,
\ldots, p_{m-1})$ and $\boldq=(q_0, q_1,\ldots, q_{m-1})$
satisfying $p_iq_i\ne 0$ for $i=0,1,\ldots, m-1$, \break one has
\begin{equation}\label{eqn6.6}
p^\ast_m = \frac{p_0p_1\cdots p_{m-1}}{p_0p_1\cdots p_{m-1}+q_0q_1\cdots
q_{m-1}} = \frac{\rho_m}{1+\rho_m}.
\end{equation}
\end{lem}

\begin{proof}
For this general case, set
\begin{equation}\label{eqn6.7}
v_m=P(A^c|S_1=1)\qquad\hbox{and}\qquad \bar v_m=P(A^c|S_1=-1).
\end{equation}
so that the expression for $p^{\ast}_m$ in (\ref{eqn2.4}) becomes
\[
p^\ast_m=p_0v_m(p_0v_m+q_0\bar v_m)^{-1}
\]
Thus (\ref{eqn6.6}) will be proved once it is established that
\begin{equation}\label{eqn6.8}
\frac{v_m}{\bar v_m}=\frac{p_1p_2\cdots p_{m-1}}{q_1q_2\cdots q_{m-1}}.
\end{equation}
By definition, $v_m$ $(\bar v_m)$ is the probability (of 'ruin')
that starting at $1$
$(-1)$ the random walk reaches $m$ $(-m)$ before it reaches $0$.
Moreover, by the modulo structure of the walk,
$\bar v_m$ is the same as the probability that starting at $m-1$, the random
walk reaches $0$ before $m$.
Thus, $v_m$, for example is the same as $\ _0f_{1m}$ in the
usual notation for these taboo probabilities; cf. Chung (1967, Sect. I.12)
where these are derived for the random walk. Direct substitution of these
exact values would then justify (\ref{eqn6.6}). Since we only require the
ratio of these two taboo probabilities, the following mapping approach
suffices, and may be of separate interest.

We first construct a $1-1$ correspondence between
the set, $\Gamma_k$, of paths that go from $1$ to $m$ without hitting $0$
and the  set, $G_k$, of paths that go from $m-1$ to $0$ without hitting $m$. Thi
s
correspondence is a simple reversal: If $\bolds_k=(s_1, s_2, \ldots, s_k,
m)$ denotes a path in $\Gamma_k$ so that $s_1=1$, $s_k=m-1$ and $1\le s_i\le
m-1$ for $1\le i\le k$, the corresponding reversed path in $G_k$ is
\[
\boldt_k(\bolds_k)=\boldt_k=(t_1, t_2, \ldots, t_k, 0) \equiv (s_k, s_{k-1},
\ldots, s_1, 0).
\]
(The reader can visualize the reversal of a path in the illustration of
Figure 1. In fact, the result becomes fairly transparent once one recognizes
the effect on paths of flipping the time axis.)
\begin{figure}
\begin{psfrags}
\includegraphics[width=6in]
{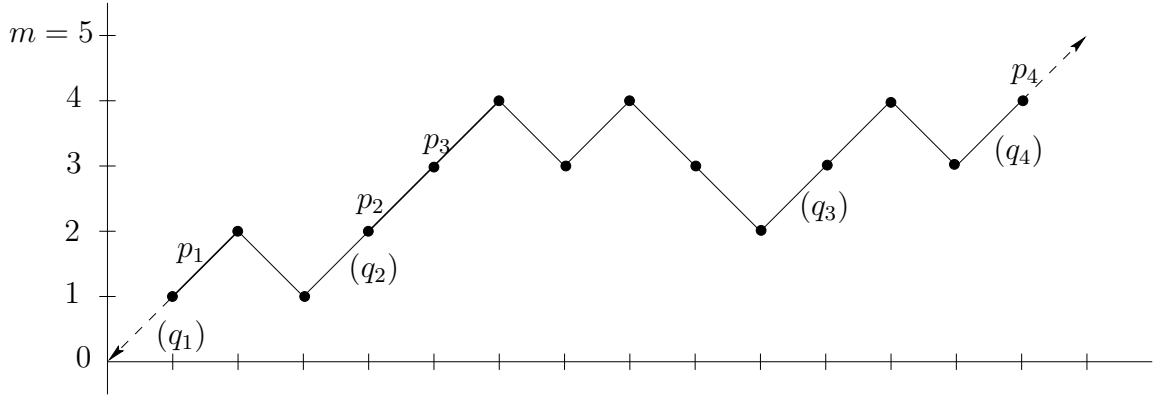}
\end{psfrags}
\caption{An illustration for $m=5$ of the correspondence between first
hitting paths from 1 to 5 and those from 4 to 0. (Probabilities in parenthesesare those for the indicated segments on the reversed path.)}
\end{figure}

For a given path $\bolds_k$, let $\nu^+_i$ ($\nu^-_i$) equal the number of
transitions from $i$ to $i+1$ ($i$ to $i-1$). Then
\begin{eqnarray}\label{eqn6.10}
P\left((S_1, S_2, \ldots, S_{k+1}) = \bolds_k|S_1=1\right) & = &
\prod^{m-1}_{i=1}p^{\nu^+_i}_i q^{\nu^-_i}_i \nonumber \\
& = & \left(\prod^{m-1}_{i=1}p_i\right)\prod^{m-2}_{j=1}
(p_jq_{j+1})^{\nu^+_j-1},
\end{eqnarray}
with the last step following since $\nu^-_1=0$, $\nu^+_{m-1}=1$ and
$\nu^-_{i+1}=\nu^+_i-1\ge 0$ for $1\le i\le m-2$. For the reversed path
$\boldt_k(\bolds_k)$, where an `up' transition of $s_j$ to $s_{j+1}$
in $\bolds_k$ becomes a `down'
transition of $s_{j+1}$ to $s_j$. Write $\widehat\nu^+_i$ and
$\widehat\nu^-_i$ for the corresponding numbers for $\boldt_k$ so that
\begin{equation}\label{eqn6.11}
P\left((S_1, S_2, \ldots, S_{k+1})=\boldt_k|S_1=m-1\right) =
\left(\prod^{m-1}_{i=1}q_i\right)\prod^{m-1}_{j=2}(q_jp_{j-1})^{\widehat\nu^-_j-
1}.
\end{equation}
But it is clear from the correspondence that $\widehat\nu^-_j=\nu^+_{j-1}$.
Thus for every $k\ge m-1$ and every path $\bolds_k\in\Gamma_k$ the ratio of
(\ref{eqn6.10}) over (\ref{eqn6.11}), namely
$p_1p_2
\cdots
p_{m-1}/q_1q_2\cdots q_{m-1}$, is constant. It now follows immediately that
(\ref{eqn6.8}) holds, thereby completing the proof.
\end{proof}

\section{Random Mixtures of Parrondo Games $\walk$}

The main question of interest for these games concerns what happens to a
player's fortune when two or more games are played in some alternating fashion.
For example, if two different games are known to be fair, can a player
create a winning game by randomly choosing between the two at each play?
Observe first of all that for $\pi\in[0,1]$, the random mixture of two
games, $\ggwalk$ and $G(,m,\bP,\bQ)$, in which at each play the former
is chosen with probability $\pi$,is also a $\mod{m}$ game, namely,
$G(m,\pi\bp+(1-\pi)\bP,\pi\bq+(1-\pi)\bQ)$. Since Lemma 3.1 characterizes
the winning or losing nature of any such game, the question of whether
the random mixture of two fair games is a winning game or not has been
theoretically answered.
By the way, the criterion in Lemma 3.1 implies that if $\bp$ and $\bq$
are interchanged in a fair game $\ggwalk$, it remains fair, whereas a losing
game would be turned into a winning game. Moreover, the nature of the
criterion is such that
it should be the exception rather than the rule for a random mixture of fair
games to remain fair. Thus at this stage, the existence of fair games whose
mixture is winning (or losing) would appear to be less paradoxical.

A couple of general questions of interest are as follows. Suppose we say
that two fair games, A and B, are
{\it mutually supportive}  if
any other game consisting of a sequence of plays of game
 A or B is not a losing
game whenever the game choices are made independently of previous outcomes.
Do mutually supportive pairs of distinct games exist?
Is it true that if a non-trivial random
mixture (in which game A is chosen independently at each stage with constant
probability) does not result in a losing game, then the two games are
mutually supportive?

In this section, we give a complete answer to the structure of random
mixtures in the special case of the
Parrondo game, $\walk$. Although this is done by rather elementary
methods, more general questions involving mixtures appear to be
quite difficult.

Consider the random mixture, $G(m, \pi p+(1-\pi)\beta, \pi p'+(1-\pi)\beta')$,
of the two Parrondo games $G(m, p, p')$ and $G(m, \beta,\beta')$, in which
the mixing probability is $\pi \in (0,1)$. Set
\begin{equation}\label{eqn3.1}
x=p/q,\ \  y=p'/q',\ \ \widehat x=\frac{\beta}{1-\beta},\ \ \widehat
y=\frac{\beta'}{1-\beta'}
\end{equation}
and
\begin{equation}\label{eqn3.2}
\overline x=\frac{\pi p+(1-\pi)\beta}{1-\{\pi p+(1-\pi)\beta\}}=
\frac{p+\lambda\beta}{q+\lambda(1-\beta)},\qquad
\overline y=\frac{p'+\lambda\beta'}{q'+\lambda(1-\beta')},
\end{equation}
where $\lambda=(1-\pi)/\pi$. Assume without loss of generality
that $\beta <p$, or equivalently, $\widehat x <x$.

The question to consider is whether the random mixture of two
losing games can be a winning game. Suppose first that the two given
games are fair. That is, by Corollary \ref{cor2.2} in the form
(\ref{eqn2.9}), our question is whether it is possible to have
\begin{equation}\label{eqn3.3}
y=x^{-m+1},\qquad \widehat y=\widehat x^{-m+1}\qquad\hbox{and}\qquad
\overline y>\overline x^{-m+1}.
\end{equation}
For simplicity, write $m-1=r$ so that $r=1,2,\cdots$. Simple algebra leads
to
\begin{equation}\label{eqn3.4}
\overline x=\frac{x(1+\widehat x)+\lambda\widehat x(1+x)}{1+\lambda+\widehat x+\lambda x},\qquad
\overline y=\frac{y(1+\widehat y)+\lambda\widehat y(1+y)}{1+\lambda+\widehat y+\lambda y}.
\end{equation}

Substitution of the first two equations of (\ref{eqn3.3}) into $\overline y$
permits the inequality $\overline y>\overline x^{-r}$ to be written after
simplification as
\begin{equation}\label{eqn3.5}
\frac{1+\lambda+\widehat x^r+\lambda x^r}{(1+\lambda)(\widehat x x)^r+\lambda\widehat x^r+x^r} >
\frac{(1+\lambda+\widehat x+\lambda x)^r}{((1+\lambda)(\widehat x x)+\lambda\widehat x+x)^r}.
\end{equation}
Clearly, this can never hold if $r=1$, (i.e $m=2$).
We assume, therefore, that $m >2$ in the remainder of this section.

If one introduces functions $f(a)=a^r$ and $g(a,b)=(1+\lambda +a+\lambda b)/((1+\lambda)ab+\lambda a+b)$, then
(\ref{eqn3.5}) involves a form of inverse composition, namely,
\[
g(f(\widehat x), f(x)) > f(g(\widehat x,x)).
\]
On the other hand, (\ref{eqn3.5}) may be written equivalently in terms
of $\pi$  as
\begin{equation}\label{eqn3.6}
\frac{1+\pi\widehat x^r+(1-\pi)x^r}{(1+\pi\widehat x+(1-\pi)x)^r}
 > \frac{1
+\pi\widehat x^{-r}+(1-\pi)x^{-r}}{1+\pi\widehat x^{-1}+(1-\pi)x^{-1})^r}.
\end{equation}
Thus, this inequality is one about norms on the simplex as may be seen as
follows: If we set $\boldu=(1,
\widehat x,x)$ and $\boldv=(1,1/\widehat x, 1/x)$,
(\ref{eqn3.5}) is equivalent to
\[
\|\boldu\|_{r,\mu}/\|\boldu\|_{1,\mu} > \|\boldv\|_{r,\mu}/\|\boldv\|_{1.\mu},
\]
where the norms are with respect to the measure $\mu$ that assigns
masses $1,\pi,1-\pi$
to the coordinates $1,2,3$, respectively.
[In this context, the special case of $\widehat x=1$, in which the first game
is the classical fair random walk, (and which is the case relevant
to the examples in
Harmer and Abbott (1999a)), is describable as a comparison between the
$r$-norms of the ray projection onto the unit simplex of the
vectors $(1,1,x)$ and $(1,1,x^{-1})$ (or equivalently, $(1,x,x)$. Moreover,
in the case of purely random
mixing $(\pi=1/2)$, the inequality is more enticing in that it may be
stated as above but for vectors $(1,1,1,x)$ and $(1,1,1,x^{-1})$ under
counting measure on the coordinates.]

Fix $\widehat x=a\ge 1$. By cross multiplying in (\ref{eqn3.6}), the
inequality is equivalent to the positivity of the polynomial
\begin{equation}\label{eqn3.7}
\begin{array}{rcl}
Q(x):&=&(1+\lambda+a^r+\lambda x^r)((1+\lambda)ax+\lambda a+x)^r \nonumber\\
&&\quad \quad -((1+\lambda)a^rx^r+\lambda a^r+x^r)(1+\lambda+a+\lambda x)^r \\
&=&(1+\lambda+a^r)((1+\lambda)ax+x+\lambda a)^r-\lambda a^r(1+\lambda +a+\lambda x)^r\\
&&\quad +x^r\{\lambda((1+\lambda)ax+x+\lambda a)^r-((1+\lambda)a^r+1)(1+\lambda +a+\lambda x)^r\} \\
&=&\D{\sum^r_{j=0}{r\choose j}\{(1+\lambda+a^r)((1+\lambda)a+1)^j(\lambda a)^{r-j}-\lambda a^r(1+\lambda +a)^{r-j}\lambda^j\} x^j}
\\
&&\D{\quad +
\sum^r_{k=0}{r\choose k}\{\lambda((1+\lambda)a+1)^k(\lambda a)^{r-k}-((1+\lambda)a^r+1)(1+\lambda +a)^{r-k}\lambda^k\}x^{r+k}.}
\end{array}
\end{equation}
Upon writing $Q(x)=\sum^{2r}_{j=0}q_jx^j$, it follows that the coefficients
are
\begin{equation}\label{eqn3.8}
q_j = \left\{ \begin{array}{ll}
{r\choose j}a^r\{(1+\lambda+a^r)(1+\lambda+a^{-1})^j\lambda^{r-j} - (1+\lambda+a)^{r-j}\lambda^{j+1}\} & \hbox{for }0\le j<r,
\\
(1+\lambda+a^r)((1+\lambda)a+1)^r - ((1+\lambda)a^r+1)(1+\lambda +a)^r & \hbox{for }j=r, \\
{r\choose j-r}a^r\{(1+\lambda+a^{-1})^{j-r}\lambda^{2r-j+1} - (1+\lambda+a^{-r})(1+\lambda+a)^{2r-j}\lambda^{j-r}\} & \hbox{for
}r<j\le 2r.\end{array}\right.
\end{equation}
Since the expressions within the parentheses in the first and third cases
are increasing in $j$, there can be at most one change of sign among the
first $r$ coefficients and at most one among the last $r$. Thus, regardless
of the sign of the middle coefficient, $q_r$, there are at most three changes
of signs in the coefficients of $Q$ with  the exact
number depending upon the signs of $q_0,q_{r-1},q_r.q_{r+1},q_{2r}$.
(One may check that $q_{2r}$ is always positive for $a>1$, while $q_0$ is
negative when $\lambda\le1$ (i.e $\pi\ge1/2i$) or when $a\le1$.)
By Descartes's rule
of signs, the number of positive roots of $Q(x)=0$ does not, therefore,
exceed 3.

It follows directly from the definition of $Q$ that $Q(a)=0$. However, one
may check that $x=a$ is in fact a double root for all positive $a$.
To see this, compute from (\ref{eqn3.7}),
\begin{equation}\label{eqn3.9}
\begin{array}{rcl}
Q'(x)&=&rx^{r-1}\lambda((1+\lambda)ax+x+\lambda a)^r\\
&&\quad + r((1+\lambda)a+1)(1+\lambda+a^r+\lambda x^r)((1+\lambda)ax+x+\lambda a)^{r-1} \\
&&\quad -((1+\lambda)a^r+1)rx^{r-1}(1+\lambda+a+\lambda x)^r\\
&&\quad -r\lambda((1+\lambda)a^rx^r+x^r+\lambda a^r)(1+\lambda+a+\lambda x)^{r-1}\\
\end{array}
\end{equation}
so that after simplification
\begin{eqnarray*}
Q'(a) &=& r(1+\lambda)^ra^{r-1}(a+1)^{r-1}\{\lambda a^r(a+1)+((1+\lambda)a+1)(1+a^r) \\
&&\quad -((1+\lambda)a^r+1)(1+a)-\lambda a(a^r+1)\}=0
\end{eqnarray*}
for any $a$. Since this implies that $x-a$ is a double root of
$Q$, it follows from Descartes's rule of signs that $Q$ has either two or three
positive roots. In either case, we need to know that the root at $x=a$
is the largest positive
root. To show this, differentiate (\ref{eqn3.9}) to obtain
\begin{eqnarray*}
r^{-1}Q''(x) &=&\lambda(r-1)x^{r-2}((1+\lambda)ax+\lambda a+x)^r + 2\lambda rx^{r-1}((1+\lambda)a+1)((1+\lambda)ax+\lambda a+x)^{r-1} \\
&&\quad +(r-1)((1+\lambda)a+1)^2(1+\lambda+a^r+\lambda x^r)((1+\lambda)ax+\lambda a+x)^{r-2} \\
&&\quad -((1+\lambda)a^r+1)\{(r-1)x^{r-2}(1+\lambda+a+\lambda x)^r+2\lambda rx^{r-1}(1+\lambda+a+\lambda x)^{r-1}\} \\
&&\quad -\lambda^2(r-1)((1+\lambda)a^rx^r+x^r+\lambda a^r)(1+\lambda+a+x)^{r-2}
\end{eqnarray*}
from which
\begin{eqnarray*}
Q''(a) &=&r(1+\lambda)^{r-1}a^{r-2}(a+1)^{r-2}\{\lambda(1+\lambda)(r-1)a^r(a+1)^2+2\lambda ra^r((1+\lambda)a+1)(a+1)\\
&&\quad +(r-1)((1+\lambda)a+1)^2(1+a^r)-(1+\lambda)(r-1)(a+1)^2((1+\lambda)a^r+1)\\
&&\quad -2\lambda ra(a+1)((1+\lambda)a^r+1)-\lambda^2(r-1)a^2(a^r+1)\}.\\
\end{eqnarray*}
By grouping the terms within the parentheses here  according to powers of $a$,
this becomes
\[
Q''(a)=(1+\lambda)^{r-1}ra^{r-2}(a+1)^{r-2}\{\lambda(r-1)(a^{r+2}-
1)+2\lambda r(a^{r+1}-a)+\lambda(r+1)(a^r-a^2)\}.
\]
Thus, for $r\ge2$ ($m\ge3$), $Q''(a)$ is positive, negative or zero
according as $a>1$, $a<1$ or $a=1$. This implies in particular that
when $a=1$,  $x=1$
is a triple root, and hence the only root by Descartes's rule of signs.
Thus, when $a=1$, $x=1$ is the only positive root, insuring that $Q(x)>0$ for
all $x>1$.
For $a>1$, the fact that $Q''(a)>0$ shows that this double
root at $x=a$ is a local minimum. Since by (\ref{eqn3.8}) the leading coefficient, $q_{2r}$,
is positive for all $\lambda$ and all $a>1$, this insures again
that $x=a$ is the
largest real root of $Q(x)=0$, thereby establishing that $Q(x)>0$ for all
$x>a$ whenever $a\ge1$. This completes the proof of

\begin{thm}\label{thm3.1}
The random mixture, $G(m, \pi p+(1-\pi)\beta, \pi p'+(1-\pi)\beta')$, of two fair games,
$G(m,\beta,\beta')$ and $G(m, p, p')$ is a winning game whenever $m\ge 3$
and $\frac{1}{2}\le\beta<p\le 1$.
\end{thm}

\begin{cor}\label{cor3.2}
There exist losing games, the random mixture of which is a winning game.
\end{cor}

\begin{proof}
By Corollary \ref{cor2.2}, the expression whose sign determines whether a
game is winning, losing or fair, is a continuous function of its variables.
It is therefore clear that for the games appearing in the statement of
Theorem \ref{thm3.1}, one may make a sufficiently small change in the
parameters $(\beta,\beta')$ and $(p,p')$ to make the associated fair games
become losing ones, while preserving the inequality that ensures that the
random mixture of the two remains a winning game.
\end{proof}

The example presented in Harmer and Abbott (1999a) may now be described as
follows. Take $m=3$, $\beta=\frac{1}{2}=\beta'$, $p=\frac{3}{4}$ and
$p'=\frac{1}{10}$. The games $G(3, \frac{1}{2}, \frac{1}{2})$ and $G(3,
\frac{3}{4}, \frac{1}{10})$ are fair by Corollary \ref{cor2.2}, so that by
Theorem \ref{thm3.1}, the mixture $G(3, \frac{5}{8}, \frac{3}{10})$
is a winning game.
Consider now the games used by these authors, $G(3, \frac{3}{4}-\e,
\frac{1}{10}-\e)$ and $G(3, \frac{1}{2}-\e, \frac{1}{2}-\e)$, and their
random mixture $G(3, \frac{5}{8}-\e, \frac{3}{10}-\e)$. It is clear that the
first two are losing games for each positive $\e<1/10$ and that
there would be some positive value
$\e_0\le1/10$ for which the mixture remains a winning game
whenever $0<\e<\e_0$, as postulated in Harmer and Abbott (1999a).

In this section we have considered the random mixing of two $\gwalk$ walks.
One is also interested in deterministic mixtures. Simulations in Harmer
and Abbott(1999a) indicate that deterministic mixtures of the two games
proposed by Parrondo turn their separate losing nature into a winning
combination. It is difficult in gneral to analyze such deterministic mixtures
since it requires computing the stationary probabilities of the
product of the associated stochastic matrices. To expand upon this, suppose
one has two distinct $\gwalk$ games called $A$ and $B$ with parameters
$a_j,0,1-a_j$ and $b_j,0,1-b_j$, respectively. By Lemma 3.2, the probabilities
$p_m^\ast(A)$ and $p_m^\ast(B)$ for the two games would equal $1/2$ (i.e.,
the games would be fair) if and only if
\begin{equation}\label{eqn3.10}
\prod^{m-1}_{j=0}\frac{a_j}{1-a_j}=1=\prod^{m-1}_{j=0}\frac{b_j}{1-b_j}.
\end{equation}
Consider now the random walk formed by alternating the transition probabilities
of these two. Then the two-step process is also a random walk, though one with
jumps of two units and with non-zero probabilities $r_i$ of zero jumps.
That is, the alternation of two $\gwalk$ games is a $\ggwalk$ game. This
$2$-step process is then reducible with two classes, the odd and the even
integers. If the walk starts in state $'0'$, for example, the corresponding
quotient of relevant  parameters is
\begin{equation}\label{eqn3.11}
\frac{(a_0b_1)(a_2b_3)\cdots (a_{m-2}b_{m-1})}{(1-a_0)(1-b_1)\cdots (1-a_{m-2})(1-b_{m-1})}.
\end{equation}
Since only half of the parameters enter here, it is clear that this ratio
may be greater or less than or equal to $1$ even when the separate games
are fair. This implies that when $m$ is even, the alternation of two
fair games may be either fair, winning or losing. Notice that even if one
imposes the natural restriction that a fair game must be fair for
all sarting states one gains nothing more since, for example,  the
condition for fairness
starting in state $'1'$ , namely,
\[
\frac{(a_1b_2)(a_3b_4)\cdots (a_{m-1}b_0)}{(1-a_1)(1-b_2)\cdots (1-a_{m-1})(
1-b_0)}=1,
\]
is equivalent under (\ref{eqn3.10}) to the expression in (\ref{eqn3.11})
being set equal to $1$.

When $m$ is odd, the alternation of fair games is fair as can be seen
by considering the two-step game as a mod $2m$ game for which fairness
requires by Lemma 3.1 that the product of (\ref{eqn3.11}) and the following
displayed quotient be equal to $1$, which follows from (\ref{eqn3.10}).
Thus the alternation of these fair games cannot result in winning ones
when $m$ is odd.

The story is different, however,  for $[AABB]$, the mixture in which
two plays of game $A$ are alternated with two plays of $B$. In view of the
previous paragraphs, this game is equivalent when $m$ is odd
to an alternating $[AB]$ game
{\it but} one in which both $A$ and $B$ are $\ggwalk$ games. For $m=3$ this
is reasonably tractable. In particular, if one of the games is the
classical simple random walk one can show that the mixture is indeed
a winning game under a natural restriction on the second game. For the
special case of $AABB$ in which $A$ and $B$ are the  fair
games $G(3, \frac{1}{2}, \frac{1}{2})$ and $G(3,
\frac{3}{4}, \frac{1}{10})$ corresponding to Parrondo's example, one can
show that the asymptotic average gain is $0.0218363>0$

\section{Direct Calculation of the Asymptotic Expected Average Gain for
a $\gwalk$ Game}

By (\ref{eqn2.8}), since $p^\ast_m$ has been evaluated,
the asymptotic average gain  (or loss) would be known once
$E(T_2-T_1)$ is computed. A closed form for this expected
inter-occurrence time is discussed below since it is of interest in its own
right for these processes. However, the asymptotic average gain,
$\lim_{n\to\infty} S_n/n$, being a limit of bounded r.v.'s, may also be
derived directly by obtaining the limit of the corresponding expectations.
We do this as follows.

Consider the game $\gwalk$. Define
\begin{eqnarray}\label{eqn4.1}
\mu^{(j)}_k:&=&E(S_{n+k}-S_n|S_n\equiv j \mod{m}) \\
&=&E(\sum^k_{i=1}E(X_i|S_0\equiv j \mod{m}),\nonumber
\end{eqnarray}
emphasizing by the notation the fact that the expectation depends only upon
the congruence class of $S_n$ modulo $m$ and not upon the actual value of $S_n$ nor of $n$. In fact, the
random walk $S_n$ is equivalent to the random walk on the circular group
of integers $\mod{m}$ where a positive move is taken to be
in the clockwise direction. Clearly,
\begin{eqnarray*}
\mu^{(0)}_{k+1} & = & p_0(1+\mu^{(1)}_k) + q_0(-1+\mu_k^{(m-1)}) \\
& = & p_0-q_0 + p_0\mu_k^{(1)} + q_0\mu_k^{(m-1)}.
\end{eqnarray*}
Similarly, for $j=1,2,\ldots, m-1$,
\begin{equation}\label{eqn4.15}
\mu^{(j)}_{k+1} = p_j-q_j + p_j\mu^{(j+1)}_k + q_j\mu^{(j-1)}_k
\end{equation}
where we equate $\mu^{(m)}_k = \mu^{(0)}_k$ and $\mu^{(-1)}_k =
\mu^{(m-1)}_k$. To express this conveniently in  matrix form, write ${\bmu}_k =
(\mu^{(0)}_k, \ldots, \mu^{r}_k)'$ and $\boldb=
(p_0-q_0, p_1-q_1, \ldots, p_r-q_r)'$ as
$m\times 1$ column vectors and set
\begin{equation}\label{eqn4.2}
\bC = \left[ \begin{array}{ccccccc}
0 & p_0 & 0 & 0 & \cdots & \cdots & q_0 \\
q_1 & 0 & p_1 & 0 & \cdots & \cdots & 0 \\
0 & q_2 & 0 & p_2 & \cdots & \cdots & 0 \\
\cdot & \cdot & \cdot & \cdot & \cdot & \cdot & \\
0 & 0 & 0 & 0 & q_{r-1} & 0 & p_{r-1} \\
p_r & 0 & 0 & 0 & 0 & q_r & 0  .
\end{array}\right]
\end{equation}
where again $r=m-1$.
Since $\bmu_1=\boldb$, it is clear from (\ref{eqn4.15}) that
\[
\bmu_k=\boldb+\bC\boldb + \bC^2\boldb + \cdots + \bC^{k-1}\boldb.
\]
which implies that
\begin{equation}\label{eqn4.3}
\lim_{n\to\infty}\bmu_n/n = \left(\lim_{n\to\infty} \frac{1}{n}
\sum^{n-1}_{i=0} \bC^i\right)\boldb.
\end{equation}

The reader should note that if $\{S_n: n \ge 0\}$ were a more general
$\mod{m}$ Markov chain, the vector $\boldb$ would be given by
\[
b_i \equiv E(X_{n+1} | S_n = i \mod{m}) = \sum_j (j-i)p_{ij}
\]
and $\bC$ would be determined by
\[
C_{ij}=P(S_{n+1}=j \mod{m} |S_n=i \mod{m}) = \sum_{k\in \bZ}p_{i,j-i+km}.
\]
That is, the transition matrix $\bC$ for the Markov chain of congruence
classes of $\{S_n\}$ is formed from the original chain's transition
matrix $\bP$ by summing over all states in the appropriate
congruence class. With these defintions, the limit of (\ref{eqn4.3})
applies to a general $\mod{m}$ Markov chain. We shall continue, however,
with the $\gwalk$ case in order to obtain explicit values.

The value of this limit depends upon the periodicity of $\bC$. Suppose
first that $m$ is {\emph odd}. In this case, $\bC$ is an irreducible
aperiodic stochastic matrix
provided only that $p_jq_j\ne 0$ for each $j$. Thus the limit exists
and is a stochastic matrix, each of whose
rows is the row vector of stationary probabilities associated with $\bC$,
$\bpi=(\pi_0, \pi_1, \ldots, \pi_{m-1})$, say. It is a known result
of G. Mihoc (cf.\ Fr\'echet
(1952), pp. 114-116) that the entries in $\bpi$ are proportional to the
diagonal cofactors of $\bI-\bC$. (See Appendix A below for this and other
results to be used below.)

Let $\gamma_{im}$ denote the $(i,i)$-th cofactor of $\bI-\bC$. These are
tractable for reasonable values of $m$. Due to the cyclic structure
underlying the matrix $\bC$ it is  necessary only to obtain
the first cofactor for each $m$. The first few values are:
\begin{eqnarray}\label{eqn6.12}
\gamma_{13}&=&1-p_1q_2,\qquad\gamma_{14}=1-p_1q_2-p_2q_3, \\
\gamma_{15}&=&1-p_1q_2-p_2q_3-p_3q_4+p_1q_2p_3q_4 \nonumber
\end{eqnarray}
and
\begin{eqnarray*}
\gamma_{16}&=&(1-p_1q_2-p_2q_3)(1-p_3q_4-p_4q_5)-p_2p_3q_3q_4\\
        &=&1-p_1q_2-p_2q_3-p_3q_4-p_4q_5+p_1q_2p_3q_4+p_1q_2p_4q_5+p_2q_3p_4q_5
\end{eqnarray*}
The remaining diagonal cofactors are then obtained for each $m$ by
successively applying the cyclic permutation of $(p_0, p_1, \ldots,
p_{m-1})$ into $(p_1, p_2, \ldots, p_{m-1}, p_0)$.
For the case of a Parrondo $\walk$ game with $m=3$, the
situation studied in Harmer and Abbott (1999a), (\ref{eqn6.12}) implies that
\[
\gamma_{13}=1-pq, \ \gamma_{23}=1-pq',\ \gamma_{33}=1-p'q.
\]

A general formula, presumably known, is possible
for these cofactors, namely,
\setlength{\arraycolsep}{1pt}
\begin{eqnarray}\label{eqn6.13}
\gamma_{1m}=1-\sum^{m-2}_{i=1}p_iq_{i+1} & + & \sum_{1\le i<j-1\le
m-3}p_iq_{i+1}p_jq_{j+1} \nonumber \\
& - & \sum_{1\le i<j-1<k-2\le m-4}p_iq_{i+1}p_jq_{j+1}p_kq_{k+1}+\cdots
\end{eqnarray}
\setlength{\arraycolsep}{5pt}
with the series continuing as long as the largest subscript does not exceed
$m-1$. Thus for $l=[(m-1)/2]$, the last term has sign $(-1)^l$ and involves
$l$ subscripts $i_1, \ldots, i_l$ satisfying
\[
1\le i_1<i_2-1<i_3-2<\cdots i_l-l+1\le m-l.
\]
As indicated by its appearance, (\ref{eqn6.13}) follows from an
inclusion-exclusion argument based on the number of pairs of adjacent
diagonal $1$'s used in the evaluation of the cofactor's determinant.
(All diagonal cofactors are of course equal for each value of
$m\ge 3$ whenever the parameters $p_j$ and $q_j$ do not depend on $j$.)

As mentioned earlier, the stationary probabilities associated with $\bC$
are proportional to these diagonal cofactors so that in our previous
notation $\pi_i=\gamma_{i+1,m}/\gamma_{\cdot m}$ where $\gamma_{\cdot m} =
\gamma_{1m} + \cdots + \gamma_{mm}$.

An early reference for the study of the general cyclical random walk on the
integers modulo $m$, the one whose transition matrix is $\bC$, is Fr\'echet
((1952), pp.\ 122--125. This is in effect a 1938 reference for this random
walk, called by Fr\'echet, ``mouvement circulaire'', since the material is
present in the 1938 first edition of his book. He works out as an example
the stationary probabilities for the case of $m=4$. He obtains $\gamma_{14}$
as $p_2p_3+q_1q_2$ which is easily seen to agree with the expression given
above in (\ref{eqn6.12}).

The asymptotic average gain given by (\ref{eqn4.3}) now follows
directly from the above for the case when $m$ is odd. It is of the form $\lambda_m(1, 1, \ldots,
1)'$ with
\begin{equation}\label{eqn4.7}
\lambda_m=\bpi_m\boldb\equiv \frac{1}{\gamma_{\cdot m}} \sum^m_{i=1}\gamma_{im}
(p_{i-1}-q_{i-1}).
\end{equation}

Consider now the case of $m$ \emph{even}, say $m=2k$ for $k\ge 2$. Then $\bC$ is the
stochastic matrix of a periodic Markov chain of period 2. By clustering the
even and odd rows and columns, it may be written in the form
\begin{equation}\label{eqn4.8}
\bC=\left[\begin{array}{cc}
0 & A \\
B & 0
\end{array}\right]
\end{equation}
in which $A$ and $B$ are $k\times k$ stochastic matrices. Consequently,
\[
\bC^2=\left(\begin{array}{cc}
AB & 0 \\
0 & BA \end{array}\right), \bC^{2s} = \left(\begin{array}{cc}
(AB)^s & 0 \\
0 & (BA)^s \end{array}\right), \bC^{2s+1} = \left(\begin{array}{cc}
0 & A(BA)^s \\
B(AB)^s & 0 \end{array}\right)
\]
in which both $AB$ and $BA$ are irreducible aperiodic recurrent stochastic
matrices. If $\bdelta$, $\brho$ represent the vectors of limiting stationary
probabilities for $AB$ and $BA$, respectively, and if $D$ and $R$ are the
matrices all of whose rows are $\bdelta$ and $(\brho)$, respectively, then
\[
\lim_{s\to\infty}\bC^{2s} = \left(\begin{array}{cc}
D & 0 \\
0 & R \end{array}\right),\qquad \lim_{s\to\infty}\bC^{2s+1} =
\left(\begin{array}{cc}
0 & R \\
D & 0 \end{array}\right)
\]
and so (\ref{eqn4.3}) becomes in the case of $m$ even,
\begin{equation}\label{eqn4.9}
\lim_{n\to\infty}\bmu_{n}/n=\frac{1}{2}\left(\begin{array}{cc}
0 & R \\
D & 0 \end{array}\right) \boldb.
\end{equation}
By the result of Mihoc, the elements of the common rows $\bdelta$ and $\brho$
of $D$ and $R$ are proportional to
the diagonal cofactors of $AB$ and $BA$,  respectively. However, as shown
in the Appendix below, the diagaonl cofactors of $\bI-\bC$ are made up of
those of $\bI_{m/2}-AB$ and $\bI_{m/2}-BA$ and that
the column sums of the latter are equal and
equal to $1/2$ of the sum of the diagonal cofactors of $\bI-\bC$; cf
(\ref{eqna6}) below. In view of (\ref{eqn4.9}) it follows that
(\ref{eqn4.7}) holds true as well when $m$ is even. We summarize this as
\begin{thm}\label{thm4.1}
For the general $\gwalk$ game, with probability one,
\begin{equation}\label{eqn4.10}
\lim_{n\to\infty} \frac{S_n}{n} \equiv\lambda_m=\bpi_m\boldb\equiv \frac{1}{\gamma_{\cdot m}} \sum^m_{i=1}\gamma_{im}
(p_{i-1}-q_{i-1}).
\end{equation}
in which  the $\gamma_i$ are the diagonal cofactors of $\bI-\bC$ and
$\gamma_{\cdot m}$ is their sum.
\end{thm}

For the special case of a $\walk$ walk,
the limit of interest in (\ref{eqn4.10}) becomes
\begin{eqnarray}\label{eqn4.4}
\lambda_m & = & \{(p'-q')\gamma_{1m} + (p-q)(\gamma_{\cdot
m}-\gamma_{1m})\}/\gamma_{\cdot m} \nonumber \\
& = & 2p-1 + 2(p'-p)\gamma_{1m}/\gamma_{\cdot m}.
\end{eqnarray}
 \ From (\ref{eqn6.12}), the first few values of $\gamma_{\cdot m}$ for
a $\walk$ walk are
\begin{eqnarray*}
\gamma_{\cdot 3} & = & 3-pq-pq'-p'q = 2 + p'p^2+q'q^2 \\
\gamma_{\cdot 4} & = & 4-4pq-2pq'-2p'q = 2(1-pq)+2(p'p^2 + q'q^2) \\
\gamma_{\cdot 5} & = & 5-9pq-3pq' - 3p'q + pq(pq+2pq'+2p'q).
\end{eqnarray*}
For Game B of Harmer and Abbott (1999a), in which $m=3$,  $p=3/4-\e$ and $p'=
1/10-\e$, one obtains
\[
\gamma_{13}=13/16-\e/2+\e^2,\qquad \gamma_{\cdot 3}=
\frac{169}{80}-\frac{\e}{5}+3\e^2,
\]
from which the limit in (\ref{eqn4.4}) becomes
\begin{equation}\label{eqn4.5}
\lambda_3 = -2\e \frac{147-24\e+240\e^2}{169-16\e+240\e^2} \cong
-1.74\e-.16 \e^2+O(\e^2)
\end{equation}

This value appears to differ from the one implied by the simulated curve for
Game B shown in Fig.\ 3 of Harmer and Abbott (1999a). The value for the curve
given there for $n=100$ is approximately
$-1.35$, whereas for $\e=.005$ and
$n=100$, the value from (\ref{eqn4.5}) is  approximately
$n\lambda_3\cong-1.74/2=-.87$. The difference is that the slope of the
simulated curve is affected by the early transient behavior; in a private
communication, Harmer and Abbott confirm the agreement with this
theoretical limit of their simulated slope
when the first 100 plays are excluded.
The analogous value for their Game A (where $p=p'=
\frac{1}{2}-\e$) is $n\lambda_3=(-2\e)n=-1$ which agrees with the curve for
Game A given in their Fig.\ 3.

For the randomized game that chooses between Games A and B
with probability $1/2$, one obtains $p=\frac{5}{8}-\e$ and $p'
= \frac{3}{10}-\e$ for which
\[
\gamma_{13} = \frac{49}{64} - \frac{\e}{4} + \e^2, \qquad \gamma_{\cdot 3} =
\frac{709}{320}-\frac{\e}{10}+3\e^2.
\]
Thus in this randomized case the asymptotic slope of $S_n/n$ is by (\ref{eqn4.4})
\begin{equation}\label{eqn4.6}
\lambda_3 = \frac{1}{4} - \frac{13\times 49}{4\times 709} - \e\left\{2 -
\frac{52\times 611}{(709)^2}\right\}+ O(\e^2) \cong .0254 - 1.9368\e + O(\e^2);
\end{equation}
the expansion used in the first step requires only that $\e < .876$.
For the parameters $n=100$ and $\e=.005$ of Fig.\ 3 of Harmer and Abbott
(1999a) the asymptotic approximation becomes $n\lambda_3\cong 2.54-.98=1.57$.
This differs from their simulated value of about 1.26, again due to early
outcome effects. The reader might note that the graphs in the insert of Fig.\ 3
seem to be closer to those of (\ref{eqn4.5}) and (\ref{eqn4.6}).

As an illustration for even $m$, consider $m=4$ for which the matrices become
\[
A=\left(\begin{array}{cc}
p' & q' \\
q & p \end{array}\right),\qquad B=\left(\begin{array}{cc}
q & p \\
p & q \end{array}\right),
\]
and
\[
AB = \left(\begin{array}{cc}
p'p+q'q & -p'p-q'q \\
2pq-1 & 1-2pq \end{array}\right).
\]
Hence $\delta_{14} = (1-2pq)(1-2pq+p'p+q'q)^{-1}$ and thus
\begin{equation}\label{eqn4.13}
\lambda_4 = \frac{(p'-p)(1-2pq)}{p'p^2+q'q^2+1-pq}+p-q =
\frac{2(p'p^3-q'q^3)}{p'p^2+q'q^2+1-pq};
\end{equation}
see also (\ref{eqn5.6}) below.

In this section, we restricted consideration to $\gwalk$ games. The approach
applies as well to $\ggwalk$ games but with the simplifying zero diagonal
of $\bC$ being replaced with the $r_j$'s.

\section{Expected Interoccurrence Times of Visits to $\bZ_m$}

Set $\tau_j=E(T_1|S_0=j)$ for $j=0,\pm 1,\ldots, \pm (m-1)$
to denote the expected time of the first visit to
$\zm$ of a $\gwalk$ walk $\{S_n\}$ starting at $j$.
In the expression (\ref{eqn2.8}) for the asymptotic average gain, the
denominator $E(T_2-T_1)$ is equal to $\tau_0$. Hence, an alternate
derivation  of the asymptotic average gain would be, in view of Lemma
(\ref{lem6.2}),
to derive $\tau_0$. This may be done by solving the recursion relations
satisfied by the $\tau_j$'s, namely,
\begin{equation}\label{eqn5.1}
\tau_j  =  p_j\tau_{j+1}+q_j\tau_{j-1}+1, \qquad \hbox{for }  j=0,\pm 1,\ldots, \pm (m-1),
\end{equation}
with boundary conditions $\tau_{-m} = \tau_m = 0$, where for negative $j$
we have $p_j=p_{j+m}$ and $q_j=q_{j+m}$ for a $\mod{m}$ walk.
The solution of (\ref{eqn5.1}) is given for example in Chung(1967, I.12.(8))
in which the reader should note that the $\rho_i$'s in this
reference are related to
the \emph{reciprocals} of those used here.

The expression that one obtains in this way is quite complicated even in the
case of $m=3$ and difficult to simplify into the more tractable expressions
that can be obtained by direct solution of (\ref{eqn5.1}) by matrix inversion.
For if $\btau:=(\tau_{m-1}, \ldots, \tau_{1}, \tau_0, \tau_{-1}, \ldots,
\tau_{-m+1})'$ is the $(2m-1)$-dimensional column vector of expected occurrence
times, $\boldone$ is the $(2m-1)$-dimensional
column vector of ones and $\bG$
denotes the $(2m-1)\times(2m-1)$ matrix of coefficients in (\ref{eqn5.1})
then the system (\ref{eqn5.1}) may be expressed as
$\btau=\bG\btau+\boldone$ whose solution, with $\bH \equiv \bI-\bG$
is expressible by
\begin{equation}\label{eqn5.3}
\btau=(\bI-\bG)^{-1}\boldone = \bH^{-1}\boldone.
\end{equation}
Thus the expected interoccurrence times of $\bZ_m$
are given as the row sums of the matrix $\bH^{-1}$.
The matrix $\bH$ whose inverse is needed is a Jacobi matrix
with $-p_i$'s below a diagonal of $1$'s  and $-q_i$'s above it, namely,
\[
\bH  = \left[\begin{array}{ccccccccc}
1 & -q_{m-1} & 0 & && & 0 & 0 & 0 \\
-p_{m-2} & 1 & -q_{m-2} &  \\
& \cdot & \cdot & \cdot & \cdot & \cdot &\\
0 & & &   1 & -q_1 & 0 & & & 0 \\
0 & & &    -p_0 & 1 & -q_0  && & 0 \\
0 & & &    0 & -p_{m-1} & 1 &  & & 0 \\
&&& \cdot & \cdot & \cdot & \cdot & \cdot & \cdot \\
0 &&&&&&-p_2 & 1 & -q_2 \\
0 &&&&&&0 & -p_1 & 1 \end{array}\right].
\]
In particular, by (\ref{eqn2.8}) the required quantity, $E(T_2-T_1)=\tau_0$,
in the computation of $p_m^\ast$, is
the sum of the middle row of $\bH^{-1}$.
Thus, if $H_{i,j}$ denotes the $\{i,j\}$-cofactor of $\bH=\bI-\bG$
and $|\bH|$  denotes the determinant of $\bH$, then $\tau_0=H_{\cdot m}/
|\bH|$ where $H_{\cdot m}=H_{1m}+\cdot \cdot \cdot+H_{2m-1,m}$.

When $m=3$, $\bH$ is a $5\times 5$ matrix whose middle cofactors
are straightforwardly shown to be
\begin{eqnarray*}
H_{13} & = & p_1p_0(1-p_1q_2),\qquad H_{23} = p_0(1-p_1q_2),\qquad H_{33}=(1-p_1q_2)^2, \\
H_{43} & = & q_0(1-p_1q_2)\quad \hbox{ and } \quad H_{53}=q_0q_2(1-p_1q_2).
\end{eqnarray*}
Hence
\[
H_{\cdot 3}=(1-p_1q_2)(3-p_1q_2-p_2q_0-p_0q_1)
\]
and
\[
|\bH|=(1-p_1q_2)(1-p_1q_2-p_2q_0-p_0q_1)=(1-p_1q_2)(p_0p_1p_2+q_0q_1q_2).
\]
Therefore, for $m=3$,
\[
\tau_0=E(T_2-T_1)=1+2/(p_0p_1p_2+q_0q_1q_2).
\]
By (\ref{eqn2.8}) and (\ref{eqn2.2}), this implies that with probability one,
\begin{equation}
\lambda_3 = \lim_{n\to\infty}\frac{S_n}{n}
 =  \frac{3(p_0p_1p_2-q_0q_1q_2)}{2+p_0p_1p_2+q_0q_1q_2}.
\end{equation}
The reader may check that this agrees with the expression  given for
$\lambda_3$ in  (\ref{eqn4.10}).

For $m=4$, $\bG$ is a $7\times 7$ matrix, and the middle column's cofactors
for the corresponding $\bH\equiv(\bI-\bG)^{-1}$ are easily computed to be
\[
\begin{array}{lll}
H_{14}=p_0p_1p_2|\bK|, & H_{24}=p_0p_1|\bK|,
& H_{34}=p_0(1-p_2q_3)|\bK|,\\[2mm]
H_{44}=|\bK|^2 , &
H_{54}=q_0(1-p_1q_2)|\bK|, & H_{64}=q_0q_3|\bK| \\[2mm]
H_{74}=q_0q_3q_2|\bK|.
\end{array}
\]
where $\bK$ is the upper left (\emph{ and} lower right) $(m-1)\times(m-1)$
corner matrix of $\bH$. This gives
\[
H_{\cdot 4} = |\bK|[3-p_0q_1-p_1q_2-p_2q_3-p_3q_0+(p_1-q_3)(p_2-q_0)]
\]
and, by expansion along the middle column, the determinant of $\bH$ is
\[
|\bH|=|\bK|[p_0p_1p_2p_3+q_0q_1q_2q_3].
\]
Therefore, after simplification,
\begin{equation}\label{eqn5.5}
\tau_0 = \frac{\bH_{\cdot 4}}{|\bH|} =
\frac{2(p_0p_1+p_2p_3+q_0q_3+q_2q_1)}{p_0p_1p_2p_3+q_0q_1q_2q_3}.
\end{equation}
So that by (\ref{eqn2.8}) and (\ref{eqn2.2})
the asymptotic slope of the random walk for $m=4$ is
\begin{equation}\label{eqn5.6}
\lambda_4=\lim_{n\to\infty}\frac{S_n}{n} = \frac{2(p_0p_1p_2p_3-q_0q_1q_2q_3)}
{p_0p_1+p_2p_3+q_0q_3+q_2q_1}
\end{equation}
with probability one. This is consistent with the result obtained by the
methods of Section 4; see (4.10).

The above discussion focuses on $\gwalk$ games rather than the more
general $\ggwalk$ games. Only minor modifications for the latter are needed.
The term $r_j\tau_j$ is added to the right hand side of the equations
(\ref{eqn5.1}). This results in a substitution of $p_j/(p_j+q_j)$ and
$q_j/(p_j+q_j)$ for the parameters of the walk, and, more significantly,
a replacement of the vector $\bI$ in the solution (\ref{eqn5.3}) by the
vector of  the reciprocals, $p_j+q_j$. A benefit of working out the
more general case would be that whenever $m$ is even, one could reduce
the problem to one of order $m/2$ by observing that the embedded walk
on $\zm$ is equivalent in its asymptotic behavior to that of the $2$-step
random walk in which the parameters would become the products, $p_0p_1,
q_0q_{m-1}$, etc. One can see this already in the example of $m=4$ above,
which the reader may compare to the case of $m=2$ for the associated
$2$-step case.

\section{A diffusion analogue of a general random walk}

Partition the real line into intervals $J_j=(j, j+1]=j+(0,1]$, for
$j=0, \pm 1, \pm 2, \cdots$. Let $\bmu=\{\mu_j:j=0, \pm 1, \cdots\}$ be
given constants.  For real $x$ set
\begin{equation}\label{eqn7.1}
\mu(x)= \sum_j\mu_j1_{ J_j}(x).
\end{equation}
Now define
 a diffusion $\{W_t:t\ge 0\}$ in terms of a
standard Brownian motion $\{B_t:t\ge 0\}$ by
\begin{equation}\label{eqn7.2}
dW_t=dB_t+\mu(W_t)dt,
\end{equation}
for $t>0$. For this process, introduce the probabilities of transition
between consecutive integers, namely,
\begin{equation}\label{eqn7.3}
p_j=p_j(\mu_j,\mu_{j-1})=P[W_\cdot\hbox{ hits }j+1\hbox{ before hitting
}j-1|W_0=j]
\end{equation}
and let $q_j=1-p_j$. Observe that $q_j(\mu_j, \mu_{j-1})=p_j(-\mu_{j-1},
-\mu_j)$ by reflection.

To obtain expressions for the $p_j$ in terms of the pertinent drift rates,
$\mu_j$ and $\mu_{j-1}$, we will use the scale function of the diffusion.
For this,
fix constants $a<b$ and define for $x\in[a,b]$ the first passage probabilities
\begin{equation}\label{eqn7.4}
u(x)=P[W_\cdot\hbox{ hits }b\hbox{ before }a|W_0=x].
\end{equation}
The backward equations for the Markov process $W_{\cdot}$ imply that
$u$ satisfies the second order differential equation $u''+2{\mu}u'=0$, the
solution of which is of the form
\[
u(x)=c\int^x_a \exp{\{-2\int^y_a\mu(z)dz\}}dy +b.
\]
The boundary conditions, $u(a)=0, u(b)=1$ then give
\begin{equation}\label{eqn7.5}
u(x) = \frac{\int^x_a \exp{\{-2\int^y_a\mu(z)dz\}}dy}{\int^b_a \exp{\{-2\int^y_a\mu(z)dz\}}dy}.
\end{equation}
Note that in the case of $\mu_j=\mu$ for every $j$, this becomes the
formula of Anderson(1960, Theorem 4.1); for $a<0<b$
\begin{equation}\label{eqn7.22}
P[B_t+\mu t\hbox{ hits }b\hbox{ before }a|B_0=0]
=\D{\frac{1-e^{2a\mu}}{1-e^{-2(b-a)\mu}}}
\end{equation}
when $\mu\ne 0$, and equals $\frac{1}{2}$ when $\mu=0$.

A scale function for the diffusion, a function, $S$ say, which satisfies $u(x)=\{S(x)-S(a)\}
/\{S(b)-S(a)\}$, may be deduced from (\ref{eqn7.5}) to be
\begin{equation}\label{eqn7.6}
S(x)=2\int^x_0\exp{\{-2\int^y_0\mu(z)\}}dy,
\end{equation}
the scalar $2$ being inserted for later simplicity.

For the step function $\mu$ considered here, the above may be integrated
out for all $x$. However, our interests here require $S$ only for integer values of $x=n$, and in this case, $S(0) = 0$ and
\begin{equation}\label{eqn7.7}
u(x) = \left\{
\begin{array}{ll}
\D{\sum^{n-1}_{k=0} r(\mu_k)\exp{\{-2\sum^k_{j=0}\mu_j\}}} & \hbox{if }n>0, \\
\D{}\\
\D{-\sum^{-1}_{k=n} r(\mu_k)\exp{\{2\sum^{-1}_{j=k+1}\mu_j\}}} & \hbox{if }n<0.
\end{array}\right.
\end{equation}

The desired transition probabilities $p_j$ follow directly now from
(\ref{eqn7.7}).
It suffices to consider $j=0$. Since $p_0 =u(0)$ when $b=1=-a$,
(\ref{eqn7.7}) implies  that
\begin{equation}\label{eqn7.8}
p_0\equiv p(\mu_0, \mu_{-1})\equiv \frac{S(0)-S(-1)}{S(1)-S(-1)}=\frac{r(\mu_{-1})}{r(\mu_0)e^{-2\mu_0}+ r(\mu_{-1})}
\end{equation}
where $r(u)=(e^{2u}-1)/u$ for $u\ne 0$ and $r(0)=2$. Note that $p(0,
0)=\frac{1}{2}$ as required for standard Brownian motion.  Using the fact
that $r(u)\exp{(-2u)}=r(-u)$, we summarize this as follows:
\begin{lem}\label{lem7.1}
For the diffusion defined by (\ref{eqn7.2}), the transition probabilities of
the embedded random walk on the integers that are defined by (\ref{eqn7.3})
are given by
\begin{equation}\label{eqn7.10}
p_j = \frac{\mu_j(e^{2\mu_{j-1}}-1)}{\mu_j(e^{2\mu_{j-1}}-1)+ \mu_{j-1}(1-
e^{-2\mu_j})} = \frac{r(\mu_{j-1})}{r(\mu_{j-1}) + r(-\mu_j)}
\end{equation}
for $j=0, \pm 1, \pm 2, \cdots$.
\end{lem}

It is clear that the recurrence or transience of this diffusion agrees with
that of the embedded random walk. By Section 6, this in turn depends upon
the quotients, $p_1p_2\cdots p_k/q_1q_2\cdots q_k$.
 \ From (\ref{eqn7.10}),
\begin{equation}\label{eqn7.11}
\frac{p_j}{q_j} = \frac{r(\mu_{j-1})}{r(\mu_j)} e^{2\mu_j}.
\end{equation}
Then for any $k\ge 1$,
\begin{equation}\label{eqn7.12}
\prod^k_{j=1}
\frac{p_j}{q_j}=\frac{r(\mu_0)}{r(\mu_k)}\exp\left\{2\sum^k_{j=1}\mu_j\right\}
\end{equation}
with a similar expression for negative indices. Substitution of these into
(\ref{eqn6.1}) and (\ref{eqn6.2}) would then determine recurrence or not.

It is of interest to point out that the $p_j$'s may be evaluated
directly from (\ref{eqn7.6}) without finding the scale function. To
see this, set $b=1=-a$ and let $x\in[-1,1]$.
By partitioning the event $[W_\cdot\hbox{ hits }1\hbox{ before }-1]$
according to hitting
$0$ or not before 1 and $-1$, the Markov property and Anderson's result
(\ref{eqn7.22}) yield
\begin{equation}\label{eqn7.13}
u(x) = \left\{
\begin{array}{ll}
\D{\frac{1-e^{-2x\mu_1}}{1-e^{-2\mu_1}}} +
\D{\frac{e^{-2x\mu_1}-e^{-2\mu_1}}{1-e^{-2\mu_1}}u(0)} & \hbox{if }x>0, \\
\D{}\\
\D{\frac{1-e^{-2(1-x)\mu_2}}{1-e^{-2\mu_2}}u(0)} & \hbox{if }x<0.
\end{array}\right.
\end{equation}
It therefore remains to derive $u(0)$.

For $\alpha\in(0,1]$ let
$v(\alpha)$ denote the value of $u(0)$ when the barriers at $\pm 1$ are
replaced by $\pm\alpha$. That is, $v(\alpha)$ is the probability of hitting
$\alpha$ before $-\alpha$ given the process starts at zero. By partitioning
the event of hitting 1 before $-1$ according to which of $\alpha$ or
$-\alpha$ is hit first, one obtains
\begin{equation}\label{eqn7.14}
u(0) \equiv v(1)  =  v(\alpha)u(\alpha) + [1-v(\alpha)]u(-\alpha) .
\end{equation}
Upon substitution of (\ref{eqn7.13}) and then solving for $v(1)$ one obtains
\begin{equation}\label{eqn7.15}
v(1)=\frac{f(\alpha)}{1/v(\alpha)+f(\alpha)-1}
\end{equation}
with
\[
f(\alpha)  =
\frac{(1-e^{-2\alpha\mu_0})(1-e^{-2\mu_{-1}})}{(e^{-2(1-\alpha)\mu_{-1}}-e^{-2\mu_{-1}})(1-e^{-2\mu_0})}.
\]
Observe that the limit of $f(\alpha)$  as $\alpha\searrow 0$ is
\[
f(0+) = \left(\frac{\mu_0}{\mu_{-1}}\right) \frac{e^{2\mu_{-1}}-1}{1-e^{-2\mu_0}}
= \frac{r(\mu_{-1})}{r(-\mu_0)}.
\]
By scaling, $p(\alpha)$ is the same as $v(1)$, \emph{but} with
$\mu_1,\mu_2$ replaced by $\alpha\mu_1, \alpha\mu_2$. Thus one
concludes that $p(0+)=\frac{1}{2}$. Substitution of these limits into the
right hand side of (\ref{eqn7.13}) leads to
\begin{equation}\label{eqn7.16}
v(1)\equiv p(\mu_0, \mu_{-1})= p_0 = \frac{r(\mu_{-1})}{r(\mu_2)+r(-\mu_0)}
\end{equation}
as desired.

For our interests here, consider the $\mod{m}$ \emph{shift diffusions} in
which $\mu_j=\mu_l$ whenever $j\equiv l\bmod{m}$. In this case,
 Lemma \ref{lem7.1} implies that
\begin{equation}\label{eqn7.16b}
\rho_m = \prod^{m-1}_{j=0} \frac{p_j}{q_j} = \exp\left\{2\sum^{m-1}_{j=0}
\mu_j\right\}
\end{equation}
so that by Lemma \ref{lem6.1}, the embedded $\mod{m}$ random walk, and hence
the $\mod{m}$ shift diffusion, is recurrent, transient toward $+\infty$ or
transient toward $-\infty$ according as
\begin{equation}\label{eqn7.17}
\sum^{m-1}_{j=0}\mu_j=0,\qquad >0\qquad \hbox{ or }\qquad <0.
\end{equation}
Then, by Lemma \ref{lem6.2}, the constant probability of ``success'' on $\z_m$
is
\begin{equation}\label{eqn7.18}
p^\ast_m=1/\left(1+\exp\left\{-2\sum^{m-1}_{j=0}\mu_j\right\}\right).
\end{equation}
Observe that the walk is fair $(p^\ast_m=\frac{1}{2})$ if and only if
$\mu_0+\mu_1+\cdots+\mu_{m-1}=0$.

If one is given the $p_j$'s, one may solve the system of equations given by
(\ref{eqn7.10}) for $j=0, 1, \ldots, m-1$ to find the shift rates $\mu_0,
\ldots, \mu_{m-1}$ for the associated $\mod{m}$ shift diffusion. For
example, for the random walk related to Game B of Harmer and Abbott (1999a)
in which $m=3$, $p_0=1/10$ and $p_1=p_2=3/4$, the drift rates are
\[
\mu_0=-.687032,\qquad\mu_1=2.748128, \qquad\mu_2=-2.06109.
\]
Note that these are proportional to $(-1, 4, -3)$. In fact, for any fair
game $G(3,p,p')$, the associated drift rates $(\mu_0,\mu_1,\mu_2)$ are equal
to $\mu_1(-q,1,-p)$ as given by

\begin{thm}\label{thm7.2}
If the transition probabilities, $p_0,p_1,p_2$, of a recurrent $\bmod{(m)}$
shift diffusion are known, the associated drift rates may be determined
uniquely as follows:
\item{i)}. If each $p_i$ equals 1/2, then each $\mu =0$;
\item{ii)}. If exactly one of the $p_i$'s, say $p_2$, is equal to $1/2$
then $(\mu_0,\mu_1,\mu_2)=(0,x,-x)$ with $x$ being the solution of
$p_0/q_0=(1-e^{-2x})/2x$;
\item{iii)}. If none of the $p_i$'s are equal to $1/2$, then
\[
(\mu_0,\mu_1,\mu_2)=(\frac{1}{2}\ln{w})(-(1-\theta),1,-\theta))
\]
in which $\theta =(1-q_1/p_1)/(1-p_0/q_0)$ and $w\equiv e^{2\mu_1}$
is the positive solution other than 1 of the equation
\begin{equation}\label{eqn7.20}
\alpha w-w^{\theta}+(1-\alpha)=0
\end{equation}
where $\alpha=(q_2/p_2)\theta$.
\end{thm}

\begin{proof}
We prove case iii) first.
Write $x=\mu_1$ and $y=\mu_2$. Set $a=p_1p_2/q_1q_2$ and $b=p_2/q_2$,
neither of which equals $1$. By (\ref{eqn7.11})
the equations to be solved are
\begin{equation}\label{eqn7.21}
a = r(x+y)/r(y) \quad \hbox{ and }\quad b =r(x)e^{2y}/r(y)=r(x)/r(-y).
\end{equation}
Observe that
\[
\frac{x+y}{y}a=\frac{e^{2(x+y)}-1}{e^{2y}-1}=1+
\frac{e^{2x}-1}{e^{2y}-1}e^{2y}=1+\frac{xb}{y},
\]
noting that the arguments of $r$ are not zero in this case.
Hence, if we set $u=x+y$, we must have $x=cu$ and $y=(1-c)u$ with
$c=(a-1)/(b-1)$. Substituion into the second equation of (\ref{eqn7.21})
gives
\[
b=(\frac{1-c}{c})\frac{e^{2cu}-1}{1-e^{2(c-1)u}}.
\]
By setting $w=e^{2x}=e^{2cu}$ this equation becomes
\[
\frac{\theta}{b}w-w^\theta +(1-\frac{\theta}{b})=0.
\]
which completes the proof of iii).

Case i) is clear. For ii), the constant $b$ above is equal to $1$. Since
$r$ is an increasing function, this means $x=-y$. The first equation then
becomes $a=r(0)/r(y)=2/r(y)$ which is equivalent to the equation given in
the statement of case ii).
\end{proof}

When $p$ is rational, the equation (\ref{eqn7.20}) of Theorem \ref{thm7.2}
becomes a polynomial. Here are two other examples:
For the fair game $G(3,2/3,1/5)$, $\theta =2/3$ and $\alpha=\theta/2$. The
equation that determines $\mu_1=\frac{1}{2}\log{w}$ is by (\ref{eqn7.20}),
$w-3w^{2/3}+2=0$. Upon setting
$z=w^{1/3}$,  the equation becomes $z^3-3z^2+2=0$, or, equivalently,
after factoring out $z=1$, $z^2-2z-2=0$. Its desired positive solution is
$z=1+\sqrt 3$ so that $\mu_1=(3/2)\log{(1+\sqrt 3)}$. This implies by
Theorem \ref{thm7.2} that the drift rates are
\[
\mu_0=-.502526,\qquad\mu_1=1.507579, \qquad\mu_2=-1.005053.
\]
For the fair Parrondo game $G(3,4/5,1/17)$, $\theta=4/5$ and $\alpha=1/5$
so that (\ref{eqn7.20}) becomes $w-5w^{4/5}+4=0$. With $z=w^{1/5}$, this
becomes, after factoring out $z=1$,
$z^4-4z^3-4z^2-4z-4=0$. The unique positive root (by Mathematica) is
$z=4.99357$ so that $\mu_1 = (5/2)\log{z}=4.020378$ so that by Theorem
\ref{thm7.2} the drift rates are
\[
\mu_0=-.804076,\qquad\mu_1=4.020378, \qquad\mu_2=-3.216302 .
\]
\begin{appendix}

\section{Results about stationary probabilities of Markov Chains}

We begin with a 1934 result of G. Mihoc that expresses stationary
probabilities of a finite state Markov chain in terms of cofactors:
See Fr\'echet (1952), pp.\ 114--6. (Mihoc's original paper was in Romanian,
and Fr\'echet elaborated upon it in his 1938 first edition of the cited
reference.)
Let $\PP$ be any $k\times k$ stochastic matrix. For any $j\in \{1, 2, \ldots,
k\}$, the following two determinants are equal, since the matrix in the
second is obtained from that in the first by replacing the $j$th column with
the sum of all columns: For $0\le s<1$,
\[
\Delta(s):=|s\bI-\PP|=\left|\begin{array}{ccccccc}
s-p_{11} & -p_{12} & -p_{1j-1} & s-1 & p_{1j+1}, & \cdots, & -p_{1k} \\
-p_{21} & s-p_{22} & & s-1 & \\
\vdots & & & \vdots & \\
-p_{k1} & & & s-1 & & & s-p_{kk}
\end{array}\right|
\]
and so
\begin{equation}\label{eqna1}
\lim_{s\nearrow 1} \frac{\Delta(s)}{s-1} = \left|\begin{array}{cccc}
1-p_{11} & -p_{12} & 1 & -p_{1k} \\
& & 1 \\
& & \vdots \\
-p_{k1} & & 1 & 1-p_{kk} \end{array}\right|.
\end{equation}
Observe that the left hand side does not depend upon $j$. Hence the right
hand side evaluated by expanding along the $j$-th column does not depend
upon $j$. That is, if $\Delta_{ij}$ denotes the $(i, j)$-th cofactor of
$\Delta(1)$,
\begin{equation}\label{eqna2}
\Delta_{\cdot j}:=\Delta_{1j}+\cdots+\Delta_{kj} = \Delta_{\cdot 1},\qquad
j=1,2,\ldots, k.
\end{equation}
On the other hand, direct evaluation of $\Delta(1)=|\bI-\PP|$ by expansion
along the $j$-th column gives
\[
\Delta(1) = \sum^k_{i=1}\Delta_{ij}(\delta_{ij}-p_{ij}).
\]
Since $\Delta(1)=0$ this shows that
\begin{equation}\label{eqna3}
\Delta_{jj}=\sum^k_{i=1}\Delta_{ij}p_{ij}.
\end{equation}
If $1$ is a \emph{simple} root of $\Delta(s)=0$, (when the corresponding
Markov chain has a single recurrent class) the derivative in (\ref{eqna1})
is non-zero so that the common sums $\Delta_{\cdot j}$ are non-zero. In this
case, (\ref{eqna3}) implies that for each $j$, $(\Delta_{1j}, \Delta_{2j},
\ldots, \Delta_{kj})/\Delta_{\cdot 1}$ is a solution in $\bfx$ of
\begin{equation}\label{eqna4}
\bfx = \bfx\PP,\qquad \sum^k_{i=1}x_i=1.
\end{equation}
Thus if $\PP$ is also such that (\ref{eqna4}) has a unique solution, which
is the case of $\PP^n$ converging as $n\to\infty$, these solutions must all
agree (with the common row elements of that limit) so that the numbers
$\Delta_{ij}/\Delta_{\cdot j}\equiv\Pi_j$ say, do not depend upon $j$.
Equivalently, the cofactors of $\bI-\PP$ form a matrix all of which columns
are equal whenever (\ref{eqna4}) has a unique solution.
(The reader will note the relationship to Cramer's rule for solving
simultaneous linear equtions.)

Even when $\PP^n$ does not converge the columns of cofactors are still all
the same as long as the corresponding Markov chain has only one recurrent
class. Here is the case of a periodic chain of period 2 which is neede
for this paper.

Suppose the stochastic matrix $\PP$ in the above discussion is a periodic
matrix $\bC$ of period 2 of the form given in (\ref{eqn4.8}), namely,
\[
\bC=\left(\begin{array}{cc}
0 & A \\
B & 0 \end{array}\right),
\]
in which $A$ is $r\times t$ and $B$ is $t\times r$ with $k=r+t$.
(The non-square nature
of $A$ and $B$ makes this slightly different from the $\bC$ of
(\ref{eqn4.7}).) Then,
\begin{equation}\label{eqna5}
\Delta(s)=|s\bI_k-\bC|=\left|\begin{array}{cc}
s\bI_r & -A \\
-B & s\bI_t\end{array}\right|.
\end{equation}

It is known (e.g., Rao (1973), p.\ 32) that determinants of this form can be
evaluated in two ways giving
\[
\Delta(s)=s|s\bI_t-s^{-1}BA|=s|s\bI_r-s^{-1}AB|.
\]
Therefore, for $u=s^2$,
\[
\Delta(\sqrt{u})=|u\bI_t-BA|=|u\bI_r-AB|
\]
and so
\[
\lim_{u\nearrow 1} \frac{|u\bI_t-BA|}{u-1} = \lim_{u\nearrow 1}\frac{|u\bI_r
- AB|}{u-1}.
\]
By (\ref{eqna1}) and (\ref{eqna2}) above, this means that the common
column sum of cofactors of $\bI_t-BA$ is equal to that of the column sums of
cofactors of $\bI_r-AB$. Moreover, since
\[
\lim_{s\nearrow 1}\frac{\Delta(s)}{s-1} = 2\lim_{u\nearrow
1}\frac{\Delta(\sqrt{u})}{u-1},
\]
each of these column sums is exactly half of the equal column sums of
cofactors of $\bI_k-\bC$.

It is possible also to show that the set of diagonal cofactors of $\bI-\bC$
is made up of the diagonal cofactors of $\bI_t-BA$ and $\bI_r-AB$. Write
$\bfalpha$ and $\bfbeta$ for the first row and column of $A$ and $B$,
respectively, so that
\[
A = \left(\begin{array}{c}
\bfalpha \\
A^\ast\end{array}\right)\hbox{ and }b=(\bfbeta, B^\ast).
\]
Then the cofactor $\Delta_{11}$ of $\bI_k-\bC$ is
\[
\Delta_{11}=\left|\begin{array}{cc}
\bI_{r-1} & A^\ast \\
B^\ast & \bI_t\end{array}\right| = |\bI_{r-1}-A^\ast B^\ast|.
\]
But since
\[
\bI_r-AB = \left(\begin{array}{cc}
1-\bfalpha\bfbeta & -\alpha B^\ast \\
-A^\ast\beta & \bI_{r-1}-A^\ast B^\ast\end{array}\right)
\]
it is clear that its first diagonal cofactor is also $|\bI_{r-1}-A^\ast
B^\ast|$. On the other hand if we use instead the partitioning
\[
A=(\bfalpha^\ast, A^{\ast\ast}), B=\left(\begin{array}{cc}
\bfbeta^\ast\\
B^{\ast\ast}\end{array}\right)
\]
in which $\bfalpha^\ast$ and $\bfbeta^\ast$ are the first column and first
row of $A$ and $B$, respectively, then
\[
\bI_k-\bC=\left(\begin{array}{cc}
\bI_r & (\bfalpha^\ast A^{\ast\ast}) \\
\left(\begin{array}{cc}
\bfbeta^\ast \\
B^{\ast\ast} \end{array}\right) & \bI_t\end{array}\right).
\]
Therefore, the $(r+1)$-th diagonal cofactor of $\bI-\bC$ is
\[
\left|\begin{array}{cc}
\bI_r & A^{\ast\ast} \\
B^{\ast\ast} & \bI_{t-1}\end{array}\right| =
|\bI_{t-1}-B^{\ast\ast}A^{\ast\ast}|.
\]
But since
\[
\bI_t-BA = \left( \begin{array}{cc}
1-\bfbeta^\ast\bfalpha^\ast & -\bfbeta^\ast A^{\ast\ast} \\
-B^{\ast\ast}\bfalpha^\ast & \bI_{t-1}-B^\ast A^\ast\end{array}\right)
\]
its first diagonal cofactor is $|\bI_{t-1}-B^\ast A^\ast|$ as well.

By cyclically permuting the first $j$ columns and rows when $1\le j\le r$,
or the $(r+1)$-th through $j$-th columns and rows when $r<j\le k$, the above
arguments prove that the first $r$ diagonal
cofactors of $\bI_k-\bC$ are those of $\bI_r-AB$ and the last $t$ of them
are the diagonal cofactors of $\bI_t-BA$.

In view of the above results, the stationary probability vectors $\bdelta$
and $\brho$ for $AB$ and $BA$, respectively, that were introduced for
(\ref{eqn4.8}) may be expressed in the notation of Section 4 as
\begin{equation}\label{eqna6}
\bdelta=\frac{2}{\gamma_{\cdot m}}(\gamma_{1m}, \ldots, \gamma_{km})\hbox{
and }\brho=\frac{2}{\gamma_{\cdot m}}(\gamma_{k+1,m}, \ldots, \gamma_{mm}).
\end{equation}
In particular, this verifies the equivalence of (\ref{eqn4.7}) and
(\ref{eqn4.9}), showing that (\ref{eqn4.9}) applies for  all $m$, whether
even or odd.

\end{appendix}

\bigskip
{\bf Acknowledgement}: The author is grateful to Derek Abbott and Greg Harmer
for their encouragement hnd for their helpful
comments on early drafts of this paper.

\references

\Ref
{Anderson, T. W. (1960). A modification of the sequential probability
ratio test to reduce the sample size. {\it Ann. Math. Statist.}
{\bf 31} 165-197.}

\Ref
{Chung, K. L. (1967). {\it Markov Chains} Springer-Verlag, New York.}

\Ref
{Durrett, R.; Kesten, H.; Lawler, G. (1991) Making money from fair games.
{\it Random walks, Brownian motion, and interacting
particle systems}, 255--267, {\it Progr. Probab.}, {\bf 28}, Birkhduser B,
Boston.}

\Ref
{Feller, W. (1968). {\it An Introduction to Probability Theory and Its
Applications}. Third Edition, J.\ Wiley and Sons, New York.}

\Ref
{Fr\'echet, M. (1952). {\it Recherches Th\'eoriques Modernes sur le Calcul
des Probabilit\'es}. Volume 2: {\it M\'ethode des Functions Arbitraires;
 Th\'eorie
des \'Ev\'enements en Chaine dans le cas d'un numbre fini d'\'etats possibles}.
 Second Edition, Gauthier-Villars, Paris.}

\Ref
{Harmer, G. P. and Abbott, D. (1999a). Parrondo's Paradox. {\it Statistical Science}
{\bf 14} 206--213.}

\Ref
{Harmer, G. P. and Abbott, D. (1999b). Parrondo's paradox: losing strategies
cooperate to win. {\it Nature,} {\bf 402} 864.}

\Ref
{Harmer, G. P., Abbott, D. and Taylor, P. (2000a). The paradox of parrondo's
games. {\it Proc. Roy. Soc. London} Ser. A. {\bf 456} 247-259.}

\Ref
{Harmer, G. P., Abbott, D., Taylor, P. G.  and Parrondo, J. M. R. (2000b)
Parrondo's paradoxical games and the discrete Brownian ratchet.
{\it Unsolved Problems of Noise and fluctuations (UPoN'99),}  A.
(Eds.D. Abbott and L. B. Kish,) American Inst. Physics, {\bf 511} 189-201.}

\Ref
{Karlin, S. and McGregor, J. (1959). Random walks. {\it Illinois J. Math.}
{\bf 3} 66-81.}

\Ref
{Rao, C. R. (1973). {\it Linear Statistical Inference and its Applications}.
Second Edition. J. Wiley and Sons, New York.}

\end{document}